\documentclass{article}
\usepackage[utf8]{inputenc}
\usepackage[english]{babel}
\usepackage{amsmath,amsthm,amsfonts, amssymb}
\usepackage{graphicx}
\usepackage{bm}
\usepackage{tikz}
\usepackage{multirow}

\graphicspath{{./pic/}}

\renewcommand{\div}{\mathop{\rm div}\nolimits}

\newcommand{\norm}[1]{\left\|#1\right\|}

\begin{document}

\title{Non-local Multi-continua Upscaling for Flows in Heterogeneous Fractured Media}
\author{
Eric T. Chung\thanks{Department of Mathematics, The Chinese University of Hong Kong, Shatin, New Territories, Hong Kong SAR, China (\texttt{tschung@math.cuhk.edu.hk}) }
\and
Yalchin Efendiev\thanks{Department of Mathematics \& Institute for Scientific Computation (ISC), Texas A\&M University,
College Station, Texas, USA (\texttt{efendiev@math.tamu.edu})}
\and
Wing Tat Leung\thanks{Department of Mathematics, Texas A\&M University, College Station, TX 77843, USA (\texttt{leungwt@math.tamu.edu})}
\and
Yating Wang\thanks{Department of Mathematics, Texas A\&M University, College Station, TX 77843, USA (\texttt{wytgloria@math.tamu.edu})}
\and
 Maria Vasilyeva\thanks{Institute for Scientific Computation, Texas A\&M University, College Station, TX, USA \& Department of Computational Technologies, North-Eastern Federal University, Yakutsk, Republic of Sakha (Yakutia), Russia
(\texttt{vasilyevadotmdotv@gmail.com})}
}

\maketitle

\begin{abstract}
In this paper, we propose a rigorous and accurate
 non-local (in the oversampled region) upscaling framework
based on some recently developed multiscale methods \cite{chung2017constraint}.
Our proposed method consists of identifying multi-continua
parameters via local basis functions and constructing non-local (in the
oversampled region)
transfer and effective properties. To achieve this, we 
significantly modify 
our recent work proposed within Generalized Multiscale Finite Element Method
(GMsFEM) in \cite{chung2017constraint} and derive appropriate local
problems in oversampled regions once we identify important modes
representing each continua. 
We use piecewise constant functions in each fracture network and in the matrix
to write an upscaled equation. Thus, the resulting upscaled
equation is of minimal size and the unknowns are average pressures
in the fractures and the matrix. 
We note that the use of non-local
upscaled model for porous media flows is not new, e.g., 
in \cite{Hamdi_Nonlocal},
the authors derive non-local approach. Our main contribution is
identifying appropriate local problems together with local
spectral modes to represent each continua.
The model problem for fractures assumes that one can identify fracture
networks.
 The resulting non-local equation (restricted to the oversampling region, which
is several times larger compared to the target coarse block)
has the same form as \cite{Hamdi_Nonlocal} with much smaller local regions.
We present numerical results, which show that the proposed
approach can provide good accuracy.

\end{abstract}

\section{Introduction}

\subsection{Flow-based upscaling methods}

Because of the level of detail in geological formations, some type of 
coarsening or upscaling is typically performed. 
In upscaling methods, media properties are upscaled and effective
properties are computed for each coarse block
\cite{bour84,dk92,cp95,dur91,cdgw03,weh02,durlofsky2003upscaling}. 
Computing
 effective properties involves solving local problems and equating
the averages of local integrated quantities. For example,
computing upscaled permeabilities in reservoir simulation
is typically based
on equating average fluxes between the local fine-grid solves
and the coarse-grid solves. This equality allows computing
the effective permeability fields.

In a more general upscaling setup, 
a multiple continua approach \cite{chung2017coupling} is
needed. In this approach, several effective properties are
computed for each coarse block in addition to 
modeling the transfer terms.
 This computation 
involves evaluating both effective properties and transfer
coefficients between different continua. The computations are
performed locally.

\subsection{Multiscale methods and their relation}

Similar to upscaling methods, many authors have recently studied
multiscale methods. In multiscale methods \cite{hw97,hkj12,jennylt03,ij04,eh09,cortinovis2014iterative,hughes95,bfhr97,hfmq98, npp08, jp05,maalqvist2014localization,egh12,chung2016adaptive,Hajibeygi11(2),hajibeygi2011loosely,chung2015generalizedwave}, instead
of computing the effective properties, one computes
multiscale basis functions. For single-phase upscaling,
a multiscale basis function for each coarse node is computed
via local solutions. These basis functions are further coupled
via a global formulation of the equation.
This approach is implemented within Multiscale
Finite Element Method (MsFEM) and other multiscale methods
\cite{hw97,hughes95,egh12,chung2016adaptive,apwy07,pwy02,paredes2017two,ch02,cgh09}.

To generalize this approach to more complex heterogeneities
and the multi-continua case,
Generalized Multiscale Finite Element Method (GMsFEM) is proposed
\cite{egh12,chung2016adaptive,chung2015generalizedwave,chung2017constraint}.
GMsFEM proposes a systematic approach to compute multiple
basis functions. This approach starts with a space of snapshots,
where one performs local spectral decomposition to compute multiscale
basis functions. Adaptivity can be used to select basis functions
in different regions. Each multiscale basis function represents
a continua as discussed in \cite{chung2017coupling} 
and there is no need for
coupling terms between these continua. The basis functions
for each continua are automatically identified.

The GMsFEM approach has been used jointly with localization ideas
in \cite{chung2017constraint}, where the authors propose Constraint Energy Minimizing
GMsFEM. In this approach, oversampling regions are used to compute
the multiscale basis functions. This construction takes
into account spectral basis functions to localize the computations.
The localization is restricted to $\log(H)$ layers and depends on
the contrast, which can be reduced using snapshot functions.
Moreover,
 it was shown that the approach converges
independent of the contrast and the convergence is
linear with respect to the coarse mesh size. More 
precisely,
the convergence is proportional to $H/\Lambda^{1/2}$,
where $\Lambda$ is the smallest eigenvalue
that the corresponding eigenfunction is not included
in the coarse space. Note that basis functions associated to
fractures correspond to very small eigenvalues.
The goal
of this paper is to modify this framework in an appropriate way
that is more suitable for flow-based upscaling and re-cast it 
as non-local upscaling.

\subsection{This paper}

To modify the multiscale approach presented in
\cite{chung2017constraint}, we first assume that one knows each separate fracture
network within a coarse-grid block. This is one of the drawbacks
of our method; however, such cases occur in many applications.
Next, we follow a general concept of spectral basis
functions and simply define constant functions in each
fracture network and the matrix. Because the fracture
has zero width, this procedure needs to be carefully formulated,
which is done in the paper.

Secondly, we solve local
problems in the oversampled region subject to the constraint
that the local problems vanish in fractures and the matrix.
This condition is imposed as a constraint to the local
problem and important for the localization. 
The local problems formulated for each continua (either fracture
network or the matrix phase)
simply minimize the local energy subject to
the constraint that the local solution ``vanishes'' in 
other continua except the one for which it is formulated for.
More precisely, for the continua $i$ in a block $K$,
we minimize the local oversampled problem such
that it is orthogonal to all continua except $i$ 
and an appropriate inner product with the continua $i$ is $1$.

It is important to note that the localization will not be possible
if we did not identify and separate each fracture network. This is due
to the fact that the effects of fractures are not localizable and
are global as it is well known.

Next, we use these local solutions to compute the upscaled equation.
Because the local calculations are done in an oversampled domain,
the transmissibilities are non-local and extend to the oversampled
region, which is $\log(H)$ layers around the target coarse block.
Our coarse-grid equations have a similar form to those
 \cite{Hamdi_Nonlocal};
however, 
we use different local problems in addition to multi-continua as well as 
localization. Moreover,
we show that one
can obtain an accurate solution independent of the contrast and the
mesh size. 
The resulting upscaled
equation is written in a discrete form as
\[
\sum_{j,n} T^{i,j}_{mn} (u_n^{(j)} - u_m^{(i)})=q_m^{(i)},
\]
where $ T^{i,j}_{mn}$ are nonlocal transmissibilities
for different continua $m$ and $n$, and $i,j$ correspond to different coarse blocks. 
{\it We note} that $ T^{i,j}_{mn}$ are defined in oversampled regions, which
are several times larger than the target coarse block. 
We investigate the non-local dependence
of these transmissibilities. We note that the proposed approach
modifies the framework developed in \cite{chung2017constraint} 
to derive the non-local
multiple continuum upscaled models.

We note that non-local upscaling is not new in porous media
\cite{cushman2002primer,kechagia2002upscaling,wen1996upscaling,efendiev2003generalized,efendiev2000modeling,wallstrom2002effective,durlofsky2007adaptive}
There have been many works related to non-local upscaling, particularly
for transport equations. However, even in elliptic equations, one
can obtain non-local upscaling results. Our proposed
method is motivated by the work of Jenny et al., \cite{Hamdi_Nonlocal},
where they derive non-local upscaled models. 
 We would also like  to note
a recent paper \cite{daniel_nonlocal}, where the authors derive non-local
upscaling for problems without high contrast. 
We remark that the upscaling for  
flows in fractured media requires multi-continua and
thus, to avoid the global upscaling, one needs to take into account
the fractures separately and localize their effects.
In all these papers, the global formulation
of the resulting macroscopic equations is the same with the main
difference related to computing upscaled quantities.
In this
regard, our approach differs from existing works in
the literature and address a general case of problems with high contrast
and multiple scales.

In the paper, we present some numerical results. In these examples, 
we compare our proposed upscaled model and the fine-grid models.
We compare both averages and downscaled quantities. Our numerical
results show that one can achieve a good accuracy with a small localization
and using several basis functions per coarse element (continua). 
More importantly, because the local functions are constants within
fractures and the matrix, our variables have physical properties and
they denote average pressures in each continua. This is very important for
practical simulations as in our previous GMsFEM framework,
one needs to extract physical parameters from the variables.

The paper is organized as follows. In the next section, 
Section \ref{sec:prelim}, we present some
preliminaries. In Section \ref{sec:method} and Section \ref{sec:parabolic}, we present 
our approach. Section \ref{sec:numresults} is devoted 
to numerical results. Finally, we present some conclusions.

\section{Preliminaries}
\label{sec:prelim}

We consider the single-phase flow equation
\begin{equation}\label{eq:flow}
-\div (\kappa(x) \nabla u) = g, \quad \text{in} \quad D
\end{equation}
subject to some boundary conditions. In our numerical examples,
we will consider the zero Neumann boundary condition $\nabla u \cdot n = 0$.
Here, $u$ is the pressure of flow, $g$ is the source term, and $\kappa(x)$ 
is a heterogeneous field with high contrast. 
We will be using a variational formulation of
(\ref{eq:flow}). To introduce it,
we denote by $V = H^1(D)$. The variational formulation 
is to find $u\in V$ such that
\begin{equation*}
\int_D \kappa \nabla u \cdot \nabla v  = \int_D gv,
\end{equation*}
where $g$ satisfies the compatibility condition $\int_D g = 0$.
For the zero Neumann boundary condition,
we will use $\int_D u = 0$ to ensure the well-posedness of the problem.

In this paper, we will mostly focus on applications 
to fractured media
and, for this reason, we also introduce some notations for fractured media.
For the fractured media, 
the domain $D$ can be divided into two sets of regions, that is 
\begin{equation}
D = D_m \bigoplus_i d_i D_{f,i}
\end{equation}
where $m$ and $f$ corresponds to matrix region and fracture regions respectively, and $d_i$ is the aperture of fracture $D_{f,i}$. The permeability in the matrix is $\kappa_m$, and the permeability in the $i$-th fracture is denoted by $\kappa_i$. We note that the permeabilities of matrix and fractures can differ by orders of magnitude. 

The solution of \eqref{eq:flow} is to find $u\in V$ such that 
\begin{equation}
a(u,v) = (g,v) \quad \forall v \in V,
\end{equation}
where $a(u,v) = \int_{D_m} \kappa \nabla u \cdot \nabla v + \sum_i \int_{D_{f,i}} \kappa_i \nabla_f u \cdot \nabla_f v$, $(g,v) = \int_D gv$.

\begin{figure}
	\begin{center}
        \includegraphics[width=0.4\textwidth]{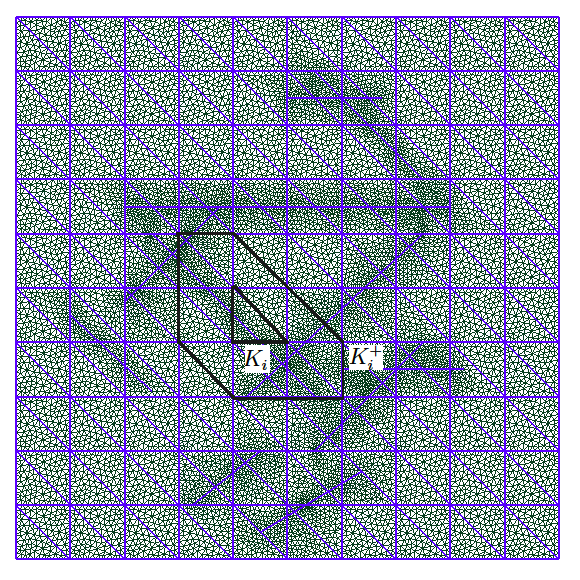}
    \end{center}    
    \caption{An illustration of coarse and fine mesh in fractured meida. $K_i$ denotes a coarse block, $K_i^1$ denotes one layer oversampling of $K_i$.} \label{fig:mesh}
\end{figure}

For the numerical approximation of problem \eqref{eq:flow}, we 
introduce the notations of fine and coarse grids.
We denote by $\mathcal{T}^H$ 
a coarse-grid partition of the domain
$D$ with mesh size $H$. 
By conducting a conforming refinement of the coarse mesh $\mathcal{T}^H$, 
we define
a fine mesh $\mathcal{T}^h$ of $D$ with mesh size $h$.
Typically, we assume that $0 < h \ll H < 1$, and 
that the fine-scale mesh $\mathcal{T}^h$
is sufficiently fine to fully resolve the small-scale information 
of the domain, and $\mathcal{T}^H$ is a coarse mesh containing 
many fine-scale features. 
We let $\{K_i | \quad i = 1, \cdots, N\}$ be the set of coarse element 
in $\mathcal{T}^H$, where $N$ is the number of coarse blocks. For each $K_i$,
the oversampled region is denoted by $K_i^+$, which is an oversampling of $K_i$ with a few layers of coarse blocks. 
An illustration of the fine and coarse meshes, as well as an oversampling region are 
shown in Figure \ref{fig:mesh}.

\section{The non-local multi-continua upscaling}

\label{sec:method}

\subsection{Multi-continua functions}

In this section, we introduce an important part of our method
that represents each continua. It will appear as 
local basis functions representing continua. In general,
these functions are automatically computed as described
in the next section, Section \ref{sec:general1}. However, for fractured
media, one can use simplified basis functions, which we discuss 
in Section \ref{sec:general2}.

\subsubsection{General spectral setup}
\label{sec:general1}

A general setup for identifying a degree of freedom for each continua
requires some spectral decomposition. Below, we briefly describe
this framework and its relation to special basis functions. We note that
this general setup shows the relevance of  special basis functions 
discussed in Section \ref{sec:general2}. Moreover,  
the general setup needs to be used when one can not identify
separate fracture networks. As we discussed that the special basis functions
do not require any basis computations provided we can identify separate fracture
networks. 

We first construct a snapshot space $V_{\text{snap}}^{i}$ for each local coarse region $\omega_i$ (a coarse neighborhood) or $K_i$ (a coarse block). 
Several choices have been developed for constructing the snapshot basis functions, including (1) the standard fine-scale basis functions, (2) harmonic basis functions which can be obtained by solving the local problems with various boundary conditons, (3) oversampling harmonic basis functions constructed in an oversampled region with standard or randomized boundary conditons, and so on. The offline space $V_{\text{off}}$ is then formed via a dimension reduction in the snapshot space using an auxiliary spectral decomposition. The offline basis obatined represent dominant modes in the snapshot space.
In the online stage, one can solve the global problem with various choices of right hand side source term or different boundary conditons. 

Here, we present a general example of constructing harmonic snapshots and spectral problem in heterogeneous media. Let $\mathcal {L}$ denotes a linear differential operator, the snapshot basis function is formed by the harmonic extension of fine-grid functions $\delta_i(x)$ which are defined on the boudary of local coarse region. That is, a snapshot basis $\psi_{\text{snap}}^{k,i}$ is the solution of 
\begin{align*}
\mathcal {L} (\psi_{\text{snap}}^{i,k}) = 0, &\quad \text{in} \quad S_i\\
\psi_{\text{snap}}^{i,k} = \delta_{i,k}, & \quad \text{on} \quad \partial S_i
\end{align*}
where $S_i \subset D$ denotes the $i$-th local region, i.e., we can take $S_i = \omega_i $ or $\omega_i^+$ or $K_i$ or $K_i^+$. Note that $\omega_i^+$ and $K_i^+$ represents the oversampled region. The fine grid functions $\delta_{i,k} = \delta_i(x_k)$ are defined for all $x_k \in \partial S$, where $\{x_k\}$ denote the fine degrees of freedom on the boundary of local coarse region $S_i$. The span of these harmonic extensions forms the local snapshot space. One can use randomized boundary conditions to reduce the computational cost associated
with snapshot calculations.

Next, we need to design a spectral problem to reduce the dimenson of local multiscale space, and the resulting space can be used as an auxiliary space for further use. Let $V(K_i)$ be the restriction of $V$ on $K_i$. Typically, we will find eigenvalues $\lambda_k^{i}$ and corresponding eigenfunction $\phi_{\text{off}}^{i,k} \in V(K_i)$ satisfying
\begin{equation}\label{eq:spectral}
a_i(\phi_{\text{off}}^{i,k}, v) = \lambda_k^{i} s_i(\phi_{\text{off}}^{i,k}, v), \quad \forall v \in V(K_i), 
\end{equation}
where bilinear operators $a_i$ and $s_i$ are defined on $V(K_i) \times V(K_i)$, and they can be symmetric non-negative definite and symmetric positive definite, respectively. For example, for the flow problem in heterogeous medium, one can choose
\begin{equation*}
a_i(u,v) = \int_{K_i} \nabla u \cdot \nabla v, \quad \quad s_i(u,v) = \int_{K_i} \tilde{\kappa} uv,
\end{equation*}
where the definition of $\tilde{\kappa}  = \sum_{j} \kappa |\nabla \chi_j|^2$ is motivated by the analysis, and $\nabla \chi_j$ denotes the multiscale partition of unity function. We arrange the eigenvalues of \eqref{eq:spectral} ascendingly, and then select the first $l_i$ eigenfunctions corresponding to the small eigenvalues to construct the offline basis functions. The span of these multiscale basis functions will form an auxiliary space, $V_{\text{aux}}^{(i)} := span\{\phi_{\text{off}}^{i,k}, \quad 1 \leq k \leq l_i\}$, where $1\leq i \leq N$ and 
$N$ is the number of coarse blocks. We note that the auxiliary space needs to be chosen appropriately, that is, all basis functions corresponding to small eigenvalues (representing the channels) have to be included in the space. 

At this point, we can construct the multiscale basis $\psi_{j,{ms}}^{(i)}$ using the auxiliary space $V_{\text{aux}}^{(i)}$ by ensuring the constraint energy minimization (see \cite{chung2017constraint} for details). We let $I_i$ be the index set containing all coarse block indices $\ell$ with $K_{\ell} \subset K_i^+$. 
To construct the required basis, we find $\psi_{j,{ms}}^{(i)}$ by solving 
\begin{equation}\label{eq:specbasis}
\begin{aligned}
\sum_{\ell\in I_i} a_{\ell}(\psi_{j,{ms}}^{(i)}, w) + \sum_{\ell \in I_i} s_{\ell}(w, \mu) &= 0, \quad \forall w \in V_0(K_i^+),\\
s_{\ell}(\psi_{j,{ms}}^{(i)}, \nu) &=  s_{\ell}(\phi_{\text{off}}^{i,j} , \nu),   \quad \forall \nu \in V_{\text{aux}}^{(\ell)}, \; \forall \ell\in I_i,
\end{aligned}
\end{equation}
where $\phi_{\text{off}}^{i,j} \in  V_{\text{aux}}^{(i)}$ is a basis in auxiliary space and $V_0(K_i^+) = H^1_0(K_i^+)$.
The contraint basis obtained form the multiscale space $V_{ms} := span \{\psi_{j,{ms}}^{(i)}, \quad 1\leq j \leq  l_i, 1 \leq i \leq N\}$, which will be used for find the multiscale solution.
Since the auxiliary space contains basis functions which capture the high-contrast features, it has been proved that the convergence of the this method is independent of the contrast and the convergence rate is in order of the coarse mesh size for appropriate oversampling size. We remark that this framework is general and can work for complex heterogeneities and multi-continuum case. For a simplified case, for example, when the fracture networks are known, we can construct some simplified basis functions with constraint energy minimization. The details are presented in the next section.

\subsubsection{Simplified basis functions representing continua}
\label{sec:general2}
In this section, we discuss the construction of simplified basis for fractured media. 
Our approach is motivated by the  Constraint Energy Minimizing
GMsFEM method proposed in \cite{chung2017constraint}, where oversampling regions are used to compute
the multiscale basis functions. We aim to construct simplified basis which has spatial decay property and can separate each contina automatically. With these simplified basis, non-local (restricted to oversampled regions)
transfer and effective properties can be constructed. 

The main idea behind this construction is to use constants within each separate
fracture network within each coarse block and a constant for each matrix.
This simplified construction of auxiliary space uses minimal degrees
of freedom in each continua. As a result, we will obtain an
upscaled equation with a minimal size. A major drawback of our
construction is that it assumes that we know the separate
fracture networks and assign a constant. Some of physical
applications can identify separate fracture networks and thus
our assumption is valid for many cases. Next, we present a detailed
description of basis construction.

We start by defining a simplified auxiliary space. 
Consider an oversampling region $K_i^+$ of coarse block $K_i$, we write $F^{(j)} = \{ f_m^{(j)}| f_m^{(j)} = D_{f,m} \cap K_j \neq \varnothing \}$ as the set of discrete fractures inside any coarse element $K_j \subset K_i^+$, and let $L_j = dim\{F^{(j)} \}$.
Let $\phi_0^{(i)},  \phi_l^{(i)} (l = 1,\cdots, L_i;  i = 1, \cdots, N)$ satisfy the following conditions  
\begin{align*}
&\int_{K_j}\phi_0^{(i)} = \delta_{ij}, \quad \int_{f_m^{(j)}}\phi_0^{(i)} = 0, \\
&\int_{K_j}\phi_l^{(i)} = 0, \quad \int_{f_m^{(j)}}\phi_l^{(i)} = \delta_{ij}\delta_{ml}.
\end{align*}
The number of fracture continuum in the coarse block $K_i$ is denoted by $L_i$. We can see that, $\phi_0^{(i)}$ has average $1$ in the matrix continua of coarse element $K_i$, and it has average $0$ in other coarse blocks $K_j \subset K_i^+$ as well as any fracture inside $K_i^+$. As for $ \phi_l^{(i)}$, it has average $1$ on the $l$-th fracture continua inside the coarse element $K_i$, and average $0$ in other fracture continua as well as the matrix continua of any coarse block $K_j \subset K_i^+$. The auxiliary space is then $V_{\text{aux}}^{(i)} = span \{ \phi_l^{(i)}, \quad 0 \leq l \leq L_i \}$. We note this definition separates the matrix and fractures, and each basis represents a continua. 

Define the subspace $ V_1(K_i^+) := \{ v \in V_0(K_i^+) | \int_K v=0,  \int_{f_m^{(j)}} v=0,   \forall K \subset K_i^+, \; \forall j,  1\leq m \leq L_j \}$. Let $G_{loc}^{(i)}: V \rightarrow V_1(K_i^+)$ be a localized operator such that 
\begin{equation*}
a(G_{loc}^{(i)}(u), v) = a(u, v), \quad \forall v \in V_1(K_i^+).
\end{equation*} 
where $V_0(K_i^+) = \{ v\in V(K_i^+) | v = 0  \text{ on }  \partial K_i^+\}$,  and $V(K_i^+)$ is the fine grid space over an oversampled region $K_i^+$.
We define $\psi_m^{(i)} := \phi_m^{(i)} - G_{loc}(\phi_m^{(i)})$, and note that the multiscale basis $\psi_m^{(i)}$ allows a spatial decay. In order to ensure constraint energy minimizing property, $\psi_m^{(i)}$ are constructed by solving the following local problem on the fine grid
\begin{equation}\label{eq:basis}
\begin{aligned}
& a(\psi_m^{(i)}, v) + \sum_{K_j \subset K_i^+} (\mu_0^{(j)} \int_{K_j}v  + \sum_{m \leq L_j} \mu_m^{(j)} \int_{f_m^{(j)}} v ) = 0, \quad \forall v\in V_0(K_i^+), \\
& \int_{K_j}\psi_m^{(i)}   = \delta_{ij} \delta_{0m}, \quad \forall K_j \subset K_i^+, \\
& \int_{f_m^{(j)}} \psi_m^{(i)}   = \delta_{ij} \delta_{nm}, \quad \forall f_m^{(j)} \in F^{(j)}, \; \forall K_j \subset K_i^+.
\end{aligned}
\end{equation}
Finally, the multiscale space for fractured media is $V_{ms} = span \{ \psi_m^{(i)}, \quad 0 \leq m \leq  L_i, 1 \leq i \leq N\}$.

\begin{figure}
	\begin{center}
        \includegraphics[width=0.3\textwidth]{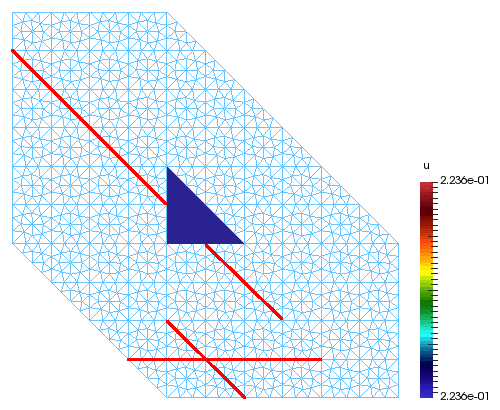}
        \includegraphics[width=0.3\textwidth]{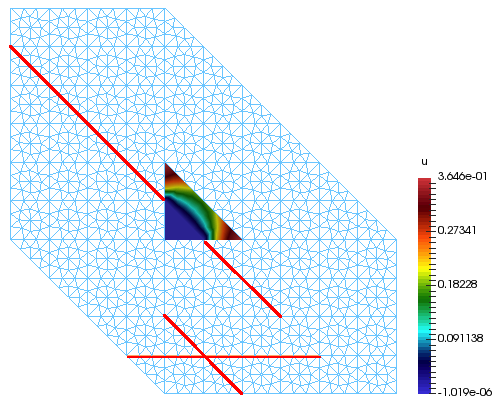}\\
        \includegraphics[width=0.3\textwidth]{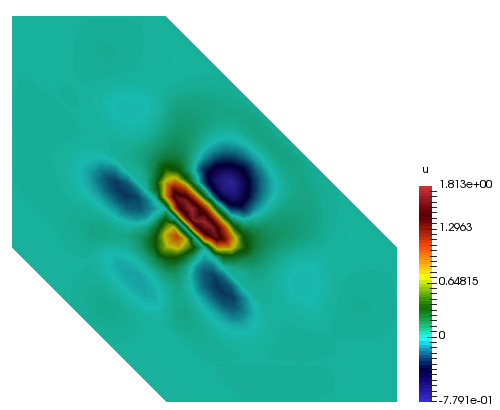}
        \includegraphics[width=0.3\textwidth]{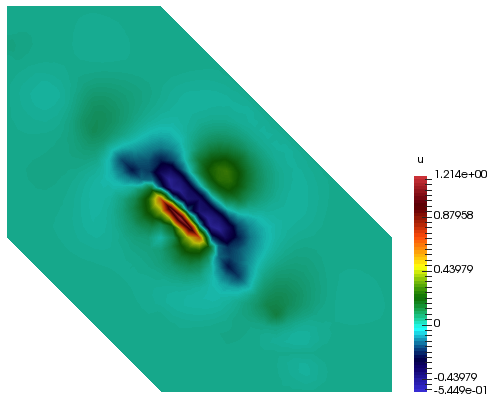}  
    \end{center}    
    \caption{Top: spectral eigenfunctions. Bottom: simplified basis.} \label{fig:spec_sim}
\end{figure}

In Figure \ref{fig:spec_sim}, we present a comparison between eigenfunctions 
constructed in Section \ref{sec:general1}, and the simplified basis constructed in this secion for a same coarse block $K$. We note that the support of the eigenbasis is in $K$, and the support of the simplified basis the oversampled region $K^+$. We can see from the top of Figure \ref{fig:spec_sim} that, the first eigenfunction is constant in the coarse block $K$, and the second eigenfunction is constant on the fracture within coarse block $K$. For the simplified basis at the bottom of Figure \ref{fig:spec_sim}, we observe that the first basis represent the matrix, and the second one represent the fracture. This indicates the relation between our simplified basis and the eigenfunctions obtained from the spectral problem. Again, we note that the simplified functions assume that one knows
separate fracture networks and uses minimal degrees of freedom to setup
a coarse system, which represent the average pressures.

\subsection{Transmissibility computations}\label{sec:trans}
\subsubsection{General spectral basis}
In the general heterogeneous case, we have constructed the constaint energy minimization basis $\{\psi_{m,{ms}}^{(i)}, \quad 1 \leq m \leq l_i, 1 \leq i \leq N\}$ as shown in \eqref{eq:specbasis}. The transmissibility matrix can be constructed by calculating
\begin{equation}\label{eq:Tloc1}
T_{mn,loc}^{(i,j)} = a(\psi_{m,{ms}}^{(i)} , \psi_{n,{ms}}^{(j)}),
\end{equation}
where $m,n$ denote the $m$- or $n$-th basis in a coarse block, $i,j$ denote the indices of coarse blocks.

\subsubsection{Simplified basis}
For the fractured media, where we assume the fractured networks are known, we constructed the simplified basis $\{\psi_m^{(i)}, \quad  0 \leq m \leq L_i, 1 \leq i \leq N \}$ by solving \eqref{eq:basis}. We define $T_{loc}$ by
\begin{equation}\label{eq:Tloc2}
T_{mn,loc}^{(i,j)} = a(\psi_m^{(i)} ,\psi_n^{(j)} ).
\end{equation}
We note that $m,n$ denotes different continua, and $i,j$ are the indices for coarse blocks. This construction shows that we can get non-local
(in the oversampled regions) transfer and effective properties for multi-continuum.

\subsection{Approximation using local multiscale basis}

Using the transmissibility defined in \eqref{eq:Tloc1} and \eqref{eq:Tloc2}, our problem is to find the approximation solution $\vec{u_T}$ such that
\begin{equation}\label{eq:approx}
\sum_n \sum_j T_{mn,loc}^{(i,j)} ([u_T]_n^{(j)} -[u_T]_m^{(i)}) = g_m^{(i)}.
\end{equation}

With a simplification of indices, we write $T_{loc}$ in the following form
\begin{equation}
\begin{bmatrix}
    t_{11} & t_{12} & \dots  & t_{1n} \\
    t_{21} & t_{22} & \dots  & t_{2n} \\
    \vdots & \vdots & \ddots & \vdots \\
    t_{n1} & t_{n2} & \dots  & t_{nn}
\end{bmatrix}
\end{equation}
where $n = \sum_{i=1}^N (1+L_i)$, and $1+L_i$ means the one matrix continua plus the number of discrete fractures in coarse block $K_i$, and $N$ is the number of coarse blocks.
The system \eqref{eq:approx} can then be expressed as in the matrix form
\begin{equation} \label{eq:finite-vol}
A_T \cdot \vec{u_T} = 
\begin{pmatrix}
   -\sum_j  t_{1j} & t_{12} & \dots  & t_{1n} \\
    t_{21} &  -\sum_j  t_{2j} & \dots  & t_{2n} \\
    \vdots & \vdots & \ddots & \vdots \\
    t_{n1} & t_{n2} & \dots  &  -\sum_j  t_{nj}
\end{pmatrix}
 \begin{pmatrix} [u_T]_0^{(1)} \\ [u_T]_1^{(1)} \\   \vdots    \\ [u_T]_N^{(L_N)} \end{pmatrix} 
 =  \begin{pmatrix}  g_0^{(1)} \\  g_1^{(1)} \\   \vdots    \\ g_N^{(L_N)} \end{pmatrix} 
\end{equation}
We remark that the summation of each row in $A_T$ is zero, which ensures the mass conservation. 

\section{Time-dependent problem}\label{sec:parabolic}

We also consider the time-dependent single-phase flow and use
spatial upscaling derived above. In particular, we consider
\begin{equation}\label{eq:unsteady1}
\frac{\partial u}{\partial t}-\div (\kappa \nabla u) = g, \quad \text{in} \quad D.
\end{equation}
The fine scale solution can be found using the standard finite element scheme, with backward Euler method for time discretization:
\begin{equation}\label{eq:unsteady2}
(\frac{u^n - u ^{n-1}}{dt} ,v) +  (\kappa \nabla u^n, \nabla v) = (g, v).
\end{equation}
In matrix form, we have
\begin{equation}
M_f u^{n} + A_f u^n = b_f + M_f u^{n-1},
\end{equation}
where $M_f$ and $A_f$ are fine scale mass and stiffness matrix respectively, $b_f$ is the right hand side vector.

For the coarse scale approximation, we will solve
\begin{equation}
M_T {u_T}^{n} + A_T {u_T}^n = b_f + M_T {u_T}^{n-1},
\end{equation}
where $A_T$ is defined in \eqref{eq:finite-vol} and $M_T$ is an approximation of coarse scale mass matrix. We note that both $A_T$ and $M_T$ are
non-local and defined for each continua. 
One can write the non-local upscaled equation as
\[
\sum_{j,n} M^{i,j}_{mn} {d\over dt} u_m^{(i)} + \sum_{j,n} T^{i,j}_{mn} (u_n^{(j)} - u_m^{(i)})=g_m^{(i)}.
\]

\section{Numerical results}
\label{sec:numresults} 
\subsection{Steady state case}

In this section, we present some representative numerical examples.
We plan to consider more realistic and complicated fracture
systems in our future works.

\begin{figure}
	\begin{center}
        \includegraphics[width=0.35\textwidth]{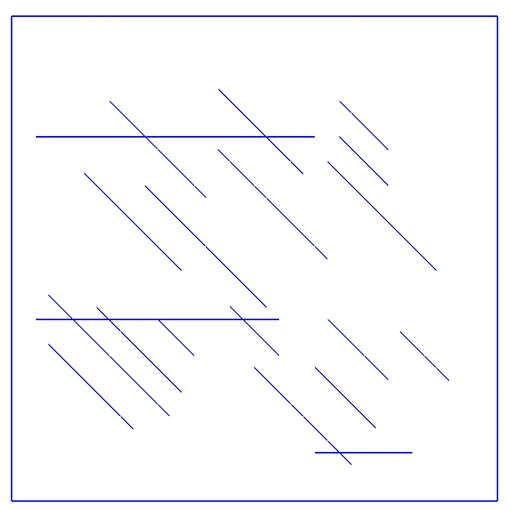}
          \includegraphics[width=0.45\textwidth]{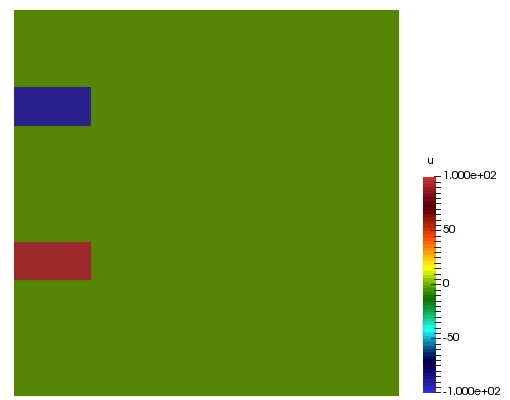}
    \end{center}    
    \caption{Left: A fractured meida. Right: source term.} \label{fig:geo_source}
\end{figure}

In this example, we use the fractured media as shown in the left of  Figure \ref{fig:geo_source}. The permeability of the matrix is $\kappa_m = 1$, and the permeability of the fractures are $\kappa_f = 10^2$. The source term in the right hand side of the equation is piecewise constant functions with $f=10^2$ for $0 \leq x \leq 0.2, 0.3\leq y \leq 0.4$, and $f=-10^2$ for $0 \leq x \leq 0.2, 0.7\leq y \leq 0.8$.

\begin{figure}
	\begin{center}
        \includegraphics[width=0.45\textwidth]{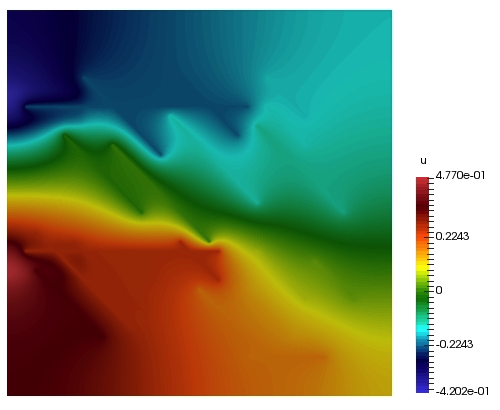}
        \includegraphics[width=0.45\textwidth]{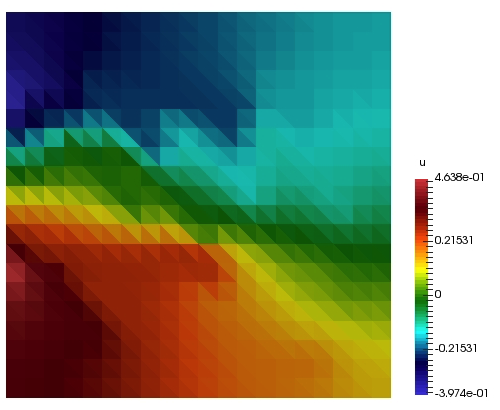}
    \end{center}    
    \caption{Left: fine scale solution. Rght: average of fine scale solution.}\label{fig:fine}
\end{figure}

The degrees of freedom for fine-scale approximation are $22642$. Let $u_f$ be the fine scale solution. We define the average of fine-scale solution $\bar{u}$ such that $\bar{u}|_{K_i }=  \int_{K_i} u, \bar{u}|_{f_m^{(i)}} = \int_{f_m^{(i)}} u$. 
These are the average pressures, which are computed with our approach.
We plot $u$ and $\bar{u}$ in Figure \ref{fig:fine}.

For the coarse scale approximation, we take the coarse mesh size as $H = 1/10, 1/20$ respectively for numerical simulations. When $H = 1/10$, the coarse scale degrees of freedom is $282$. When $H = 1/20$, the coarse scale degrees of freedom is $927$.

\begin{figure}
	\begin{center}
        \includegraphics[width=0.2\textwidth]{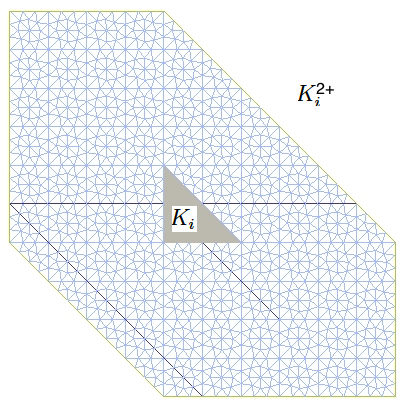}
        \includegraphics[width=0.25\textwidth]{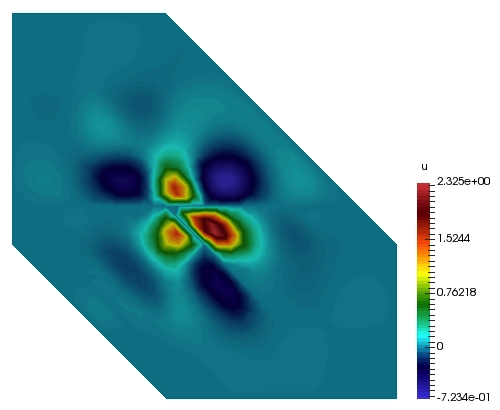}
        \includegraphics[width=0.25\textwidth]{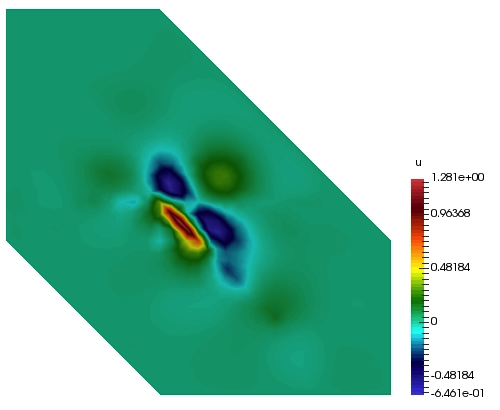}
    \end{center}    
    \caption{From left to right: A coarse block $K_i$ with two oversampling layers $K_i^{2+}$.  Local solution w.r.t matrix. Local solution w.r.t. the fracture.}\label{fig:Basis2}
\end{figure}

\begin{figure}
	\begin{center}
        \includegraphics[width=0.2\textwidth]{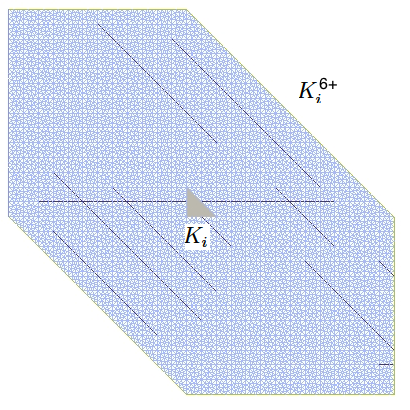}
        \includegraphics[width=0.25\textwidth]{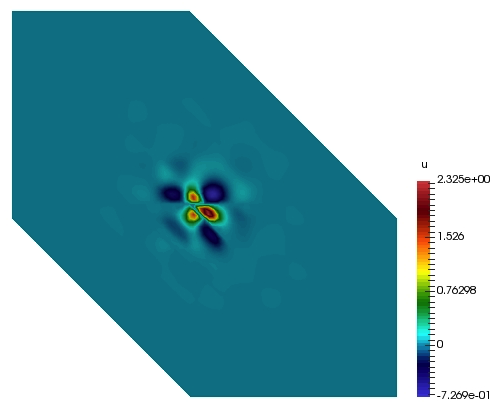}
        \includegraphics[width=0.25\textwidth]{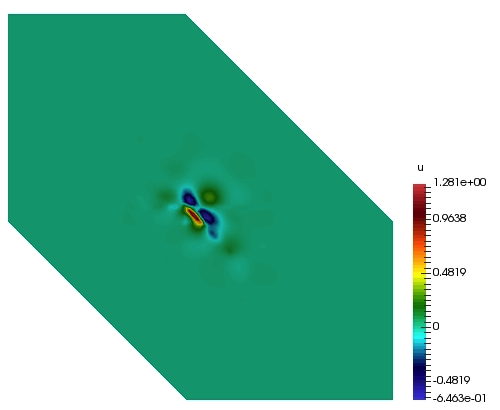}
    \end{center}    
    \caption{From left to right: A coarse blcok $K_i$ with six oversampling layers $K_i^{6+}$.  Local solution w.r.t matrix. Local solution w.r.t. the fracture.}\label{fig:Basis6}
\end{figure}

First, we present the local solutions constructed. We take the example 
for the coarse mesh with $H= 1/20$ and consider a coarse block $K_i$ 
with one discrete fracture in it. The two local solutions
 are shown in Figure \ref{fig:Basis2} when we use two oversampling layers, 
and in Figure \ref{fig:Basis6} when 
we use six oversampling layers, respectively. From the figures, 
we notice that with two layers of oversampling, the local distribution has a decay property, and the local solution almost vanish outside the oversampling 
region with six layers. This indicates that one can localize the effects.

\begin{figure}
	\begin{center}
        \includegraphics[width=0.45\textwidth]{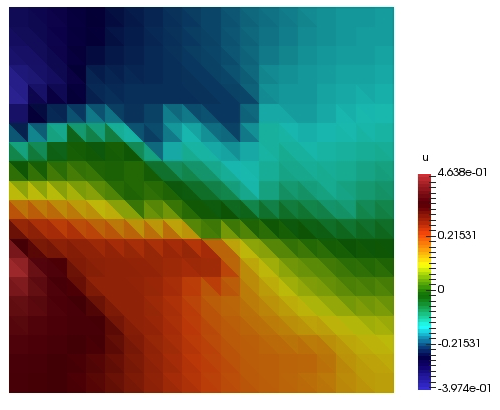}
        \includegraphics[width=0.45\textwidth]{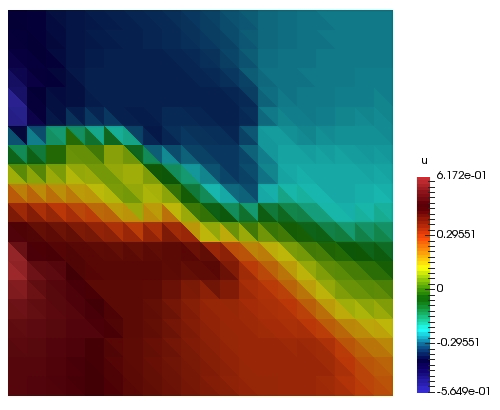}
         \includegraphics[width=0.45\textwidth]{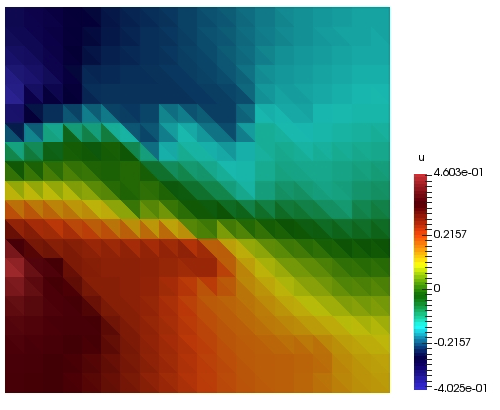}
          \includegraphics[width=0.45\textwidth]{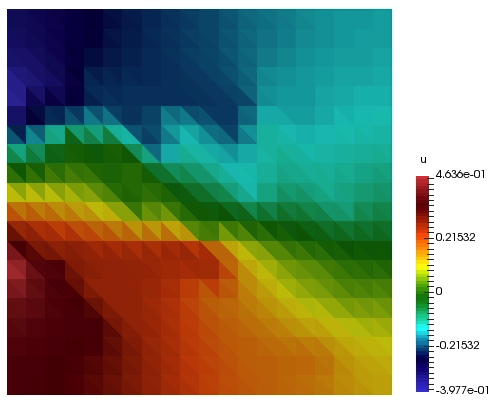}
    \end{center}    
    \caption{Corse mesh size $1/20$. Upper left: coarse scale solution using global domain for local computations. Upper right: coarse scale solution with oversampling size 1. Lower left: coarse scale solution with oversampling size 2. Lower right: coarse scale solution with oversampling size 4.}
\label{fig:ms}
\end{figure}

We present the error $\norm{u_T-\bar{u}}_{L_2}$ in Tables \ref{err-10} and \ref{err-20} for $H = 1/10$ and $H = 1/20$, respectively. From the numerial results, we observe a good convergence comparing $u_T$ with the averaged fine-scale solution. It can be seen that, with 2 layers of oversampling, the error $\norm{u_T - \bar{u}}_{L_2}$ is $1.46\%$ for $H = 1/10$. For the case $H = 1/20$, 4 layers of oversampling gives an error of $0.008\%$. This indicates the upscaled equation in our modified method can use small local regions. We plot the upscaled solutions using different size of oversampling region and compare them with the averaged fine-scale solution for the case $H=1/20$. The results are presented in Figure \ref{fig:ms}. It shows that we can obtain very good accuracy with 4 layers of oversampling.  

\begin{table}[!htb]
\centering
  \begin{tabular}{ |c  | c |}
    \hline
$Oversampling$ &$\norm{u_T - \bar{u}}_{L_2}$(\%)\\  \hline 
1 &21.97	\\  \hline
2 &1.46 \\  \hline
3 &0.015 \\  \hline
4 &0.0008\\  \hline
global &4.57e-10\\  \hline
  \end{tabular}
\caption{Coarse mesh size $1/10$. Upscaling errors when oversampling with 1,2,3,4 layers of coarse blocks. Last row shows the error when using global domain for the local computations.}\label{err-10}
\end{table}

\begin{table}[!htb]
\centering
  \begin{tabular}{ |c | c |}
    \hline
$Oversampling$ &$\norm{u_T - \bar{u}}_{L_2}$(\%)\\  \hline 
1  &40.63\\  \hline
2  &22.95\\  \hline
4  &0.008\\  \hline
6  &0.0007\\  \hline
global &6e-5\\  \hline
  \end{tabular}
\caption{Coarse mesh size $1/20$. Upscaling errors when oversampling with 1,2,4,6 layers of coarse blocks. Last row shows the error when using global domain for local computations. }\label{err-20}
\end{table}

\begin{figure}
	\begin{center}
        \includegraphics[width=0.45\textwidth]{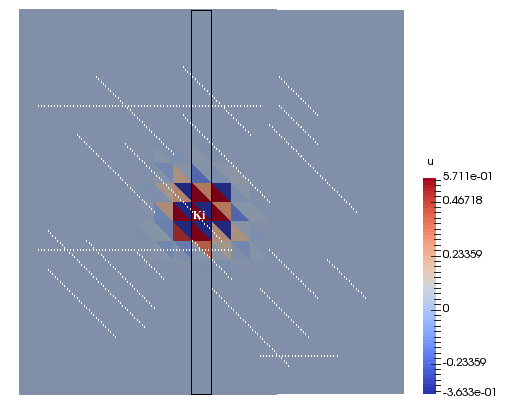}
        \includegraphics[width=0.45\textwidth]{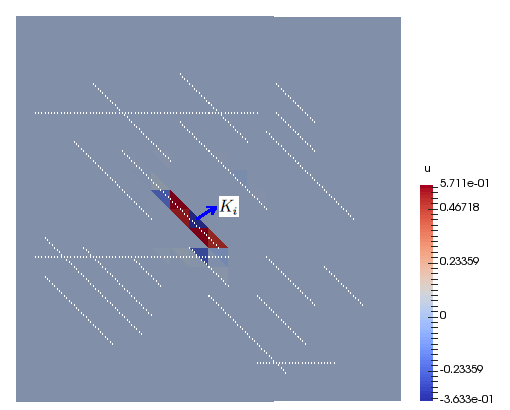}
    \end{center}    
    \caption{Using global domain for local computation. Left: Transmissibility between an element $K_i$ and neighboring elements for matrix. Right: Transmissibility between $K_i$ and neighboring elements for fractures. The dotted white lines denotes fractures in the domain.}\label{fig:over-all}
\end{figure}

\begin{figure}
	\begin{center}
        \includegraphics[width=0.55\textwidth]{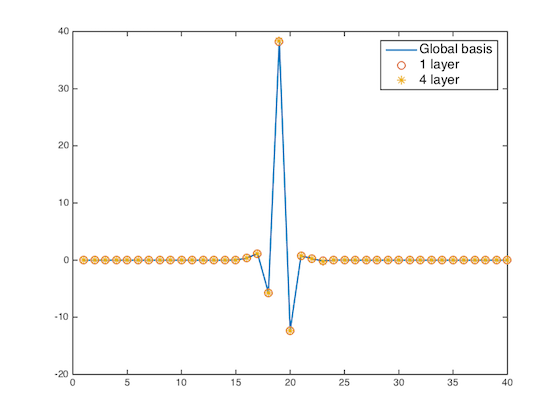}
         \includegraphics[width=0.55\textwidth]{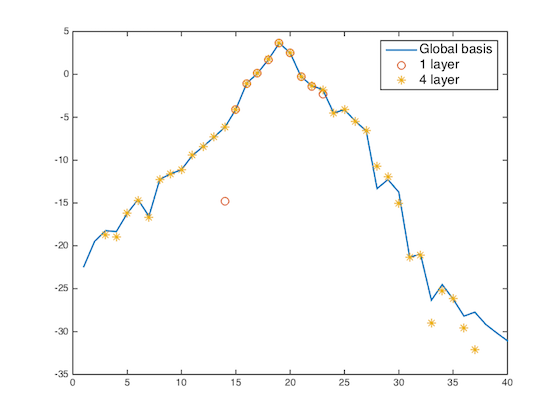}
    \end{center}    
    \caption{Transmissibility between an element $K_i$ and elements in the cross section marked by the black box in Figure \ref{fig:over-all} (left). Top: plot of exact values. Bottom: log plot of absolute values. }\label{fig:trans-1d}
\end{figure}

\begin{figure}
	\begin{center}
        \includegraphics[width=0.45\textwidth]{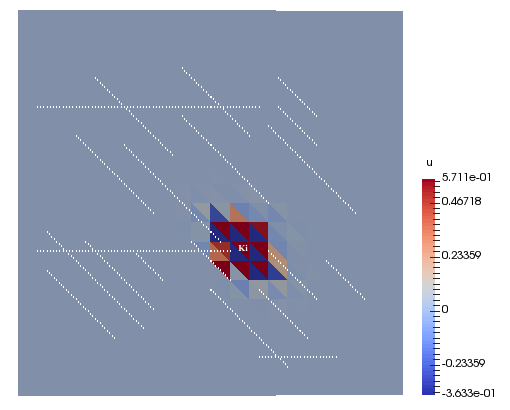}
        \includegraphics[width=0.45\textwidth]{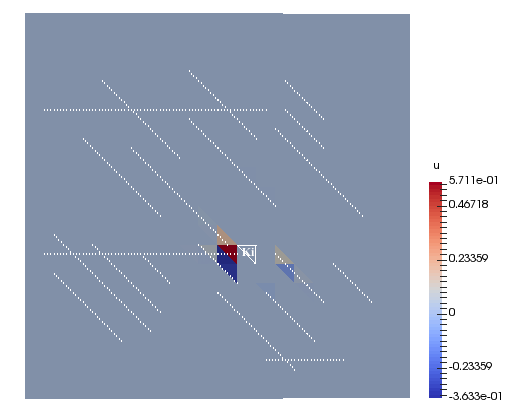}
    \end{center}    
    \caption{Using the global domain for local computations. Left: Transmissibility between an element $K_i$ and neighboring elements for matrix. Right: Transmissibility between $K_i$ and neighboring elements for fractures. The dotted white lines denotes fractures in the domain.}
\label{fig:over-all2}
\end{figure}

Next, we present numerical results for the transmissibility matrix $T_{loc}$. 
We display $T_{mn,loc}^{(i,j)}$ for two different coarse blocks $K_i$ in Figure 
\ref{fig:over-all} and Figure \ref{fig:over-all2} in the global domain. 
In the left of these two figures, we plot the transmissibility between
 the element $K_i$ and its neighboring elements for the matrix continua. 
In the right of the figures, we show the transmissibility for the 
fracture continua between $K_i$ and neighboring elements.  We notice that 
the regions of influence are almost within the 4 layers of overampling region,
which is in accordance with our numerical results. 
Figure \ref{fig:trans-1d} shows a one dimensional plot of 
transmissibility between an element $K_i$ and other coarse elements in the 
region marked by the black box in Figure \ref{fig:over-all} (left). 
We remark that the transmissibility matrix computed using global 
domain is exact, which can be used as a reference. From the subplot 
on the top of Figure \ref{fig:trans-1d}, we note that the 
transmissibility between $K_i$ and other coarse blocks along the 
slab decays fast. It can be seen from the lognormal plot on the 
bottom of Figure \ref{fig:trans-1d}, that with one layer of 
oversampling, the values in the transmissibility matrix has 
large errors, however, with four layers of oversampling, the 
transmissibility is quite accurate.

\subsection{Time-dependent case}

In this example, we present numerical results for the time-dependent case. 
The source term and geometry are the same as in the previous section, shown in Figure \ref{fig:geo_source}. 
The simulation runs for a total time of $T = 1.0$, we present the 
error between the coarse scale solution and average fine-scale solution 
at time instances $t = 0.1, 0.5$ and $1.0$. The results show good accuracy 
of the proposed method.
\begin{table}[!htb]
\centering
  \begin{tabular}{ |c | c | c | c | }
    \hline
Oversampling & $t = 0.1$ & $t = 0.5$ & $t = 1.0$ \\  \hline 
1 	&40.055 & 61.027 & 62.741	\\  \hline
2 	&1.065 & 1.035 &1.035    \\  \hline
4 	& 0.076 & 0.009  & 0.008   \\  \hline
6 	& 0.072 & 0.002  & 0.0007  \\  \hline
global  &0.072 & 0.001 &0.0006  \\  \hline
  \end{tabular}
\caption{ Error $\norm{u_T - \bar{u}}_{L_2}$(\%). Coarse mesh size $1/20$. Upscaling errors when oversampling with 1,2,4,6 layers of coarse blocks. Last row shows the error when using global domain for local computation. }\label{para-err-20}
\end{table} 

\begin{figure}
	\begin{center}
        \includegraphics[width=1.0\textwidth]{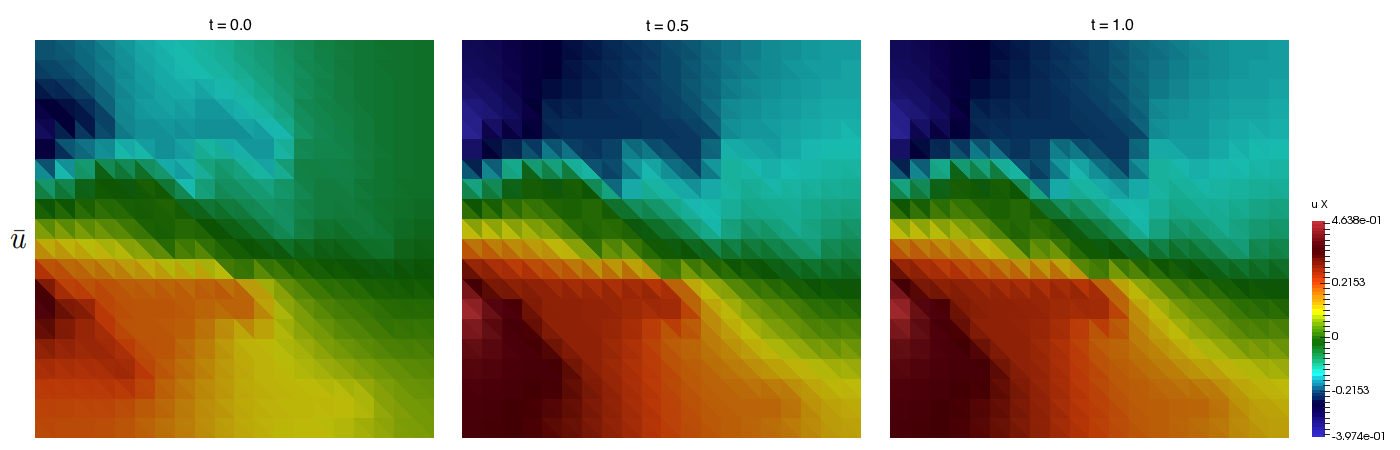}
        \includegraphics[width=1.0\textwidth]{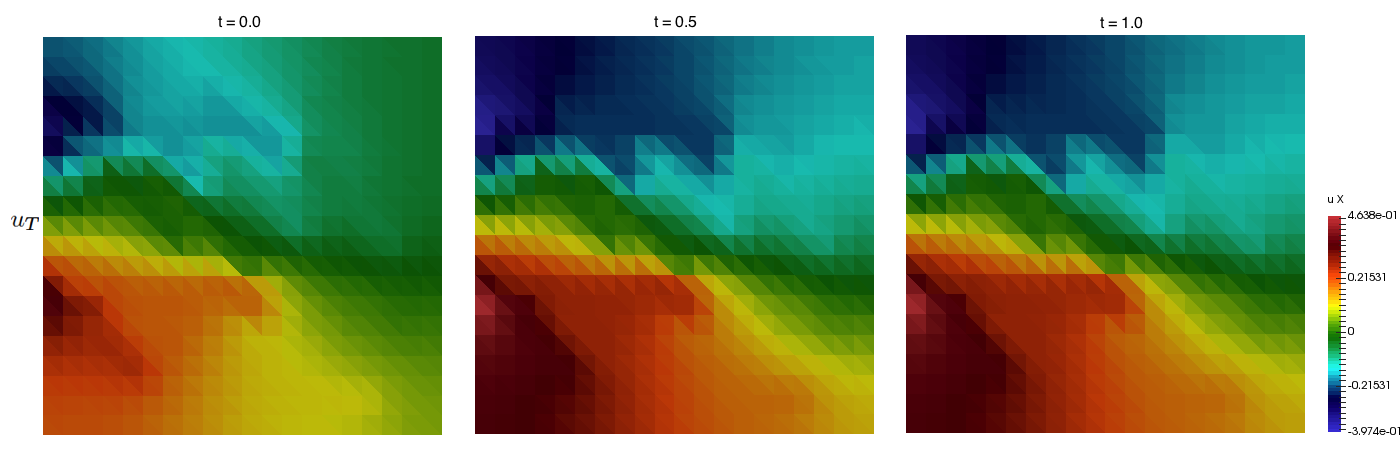}
    \end{center}    
    \caption{Oversampling size 4. Top: Average of fine scale solution at $t = 0, 0.5, 1.0$. Bottom: Coarse scale solution at $t = 0, 0.5, 1.0$.}
\label{fig:over-all22}
\end{figure}

\section{Conclusions}

We propose a non-local upscaling framework
based on some recently developed multiscale methods \cite{chung2017constraint}.
The approach uses constraints local solutions to compute the non-local
effective transmissibilities for each continua restricted to
the oversampled regions, which are several times larger than the target
coarse block.
The continua are defined by choosing piecewise constant
functions for each fracture network and the matrix.
Because of our choices of the variables,
the resulting system defined average pressures for each continua,
which is important for applications. 
The resulting nonlocal upscaled equation is defined in a small
neighborhood of the coarse block.
The local problems for the computation of effective properties
are formulated for each continua such that it is orthogonal to other
continua defined by piecewise constant functions.
We note that the use of non-local
upscaled model for porous media flows is not new, 
e.g., in \cite{Hamdi_Nonlocal},
the authors derive non-local approach.
One can use these ideas to model non-local upscaled quantities
analytically, which we will pursue in our future works.
We present numerical results, which show that the proposed
approach can provide a good accuracy. We compare average pressures
and study the decay property of local upscaled quantities.
Our numerical examples are simplistic and we plan to consider
more general fracture systems in our future works.

\section*{Acknowledgement}

The research of Eric Chung is partially supported by the Hong Kong RGC General Research Fund (Project 14317516)
and the CUHK Direct Grant for Research 2016-17.
MV's work is partially supported by Mega-grant of the Russian Federation Government (N 14.Y26.31.0013).
YE would like to thank Lou Durlofsky for the discussions related to fracture-matrix upscaling during IPAM meeting.

\bibliographystyle{plain}
\bibliography{references,references1,references_outline}

\end{document}